\documentstyle[11pt,leqno]{article}
\newtheorem{theorem}{{\sc Theorem}}
\newcommand{\bt}{\begin{theorem}}
\newcommand{\et}{\end{theorem}}
\newcommand{\newsection}[1]{\setcounter{equation}{0} \setcounter{theorem}{0}
\section{#1}}

\newcommand{\NI}{\noindent}
\newcommand{\bea}{\begin{eqnarray}}
\newcommand{\eea}{\end{eqnarray}}

\def \spec#1 {\mathop{#1}}

\def \b #1 {\bf #1}

\newcommand {\CC}{\centerline}

\newcommand{\clf}{{\cal F}}

\newcommand{\ity}{\infty}
\newcommand{\raro}{\rightarrow}

\newcommand{\vsp}{\vskip 1em}

\newcommand{\be}{\begin{equation}}
\newcommand{\ee}{\end{equation}}
\newcommand{\ben}{\begin{eqnarray*}}
\newcommand{\een}{\end{eqnarray*}}

\begin{document}
\sloppy
\CC {\Large{\bf  Nonparametric Estimation   of Trend for }}
\CC {\Large {\bf Stochastic Differential Equations Driven by}}
\CC{\Large{\bf Multiplicative Stochastic Volatility }}

\vsp \CC {B.L.S. PRAKASA RAO}
\CC{CR RAO Advanced Institute of Mathematics, Statistics }
\CC{and Computer Science, Hyderabad, India}
\vsp 
\NI{\bf Abstract:} We discuss nonparametric estimation of the trend coefficient in models
governed by a stochastic differential equation driven by a multiplicative stochastic volatility.
\vsp 
\NI{\bf Keywords and phrases :} Trend coefficient; Nonparametric estimation; Kernel method; Multiplicative stochastic volatility;
Brownian motion. 
\vsp
\newsection{Introduction} Khlifa et al. (2016) studied a class of stochastic differential equations of the form 
$$dX_t=a(t,X_t)dt+\sigma_1(t,X_t)\sigma_2(t,Y_t)dW_t, X_0=x_o\in R, 0\leq t \leq T$$
where $\{W_t, t\geq 0\}$ is a standard Wiener process and $\{Y_t, t \geq 0\}$ is a stochastic process adapted to the filter generated by the Wiener process $W.$ They proved the following results concerning the existence and uniqueness of the solution for such SDEs.
\vsp
\NI{\bf Theorem 1.1:} {\it Let $Y$ be an adapted continuous process, $a(t,x), \sigma_1(t,x)$ and $\sigma_2(t,x)$ be continuous functions with respect to $x$ and $y$ in $R$ for $t\in [0,T], \sigma_2 (t,x)$ be bounded, and
\be
|\sigma_1(t,x)|^2+ |a(t,x)|^2\leq K(1+|x|^2), 0\leq t , x \in R 
\ee
for some constant $K>0.$ Then the equation (1.1) has a weak solution.}
\vsp
\NI{\bf Theorem 1.2:} {\it Let $Y$ be an adapted continuous process and $a(t,x), \sigma_1(t,x)$ and $\sigma_2(t,x)$ be  functions such that 
the following conditions hold.\\
(i) there exists a positive increasing function $\rho(u), u \in (0,\infty),$ satisfying $\rho(0)=0$ such that
$$|\sigma_1(t,x)-\sigma_1(t,y)|\leq \rho(|x-y|), t \geq 0, x,y \in R; \int_0^\ity \rho^{-1}(u)du=+\infty;$$
and\\
(ii) there exists a positive increasing concave function $k(u), u\in (0,\ity),$ satisfying $k(0)=0$ such that 
$$|a(t,x)-a(t,y)|\leq k(|x-y|), t \geq 0, x,y \in R; \int_0^\ity k^{-2}(u)du=+\infty.$$ 
Then the equation (1.1) has a pathwise unique solution and hence has a unique strong solution.}
\vsp
\NI{\bf Theorem 1.3:} {\it Let $Y$ be an adapted continuous process and $a(t,x), \sigma_1(t,x)$ and $\sigma_2(t,x)$ be  functions such that the following conditions hold:\\
(i) there exists a constant $K>0$ such that , for all $t\geq 0,$ and $x \in R,$
$$|\sigma_1(t,x)|^2+|a(t,x)|^2\leq K(1+|x|^2);$$
(ii) for any integer $N\geq 1,$ there exists a constant $ K_N>0$ such that, for all $t \geq 0$ and for all $x,$ satisfying $|x|\leq N$ and $|y|\leq N,$
$$|a(t,x)-a(t,y)|+|\sigma_1(t,x)-\sigma_1(t,y)|\leq K_N|x-y|;$$
and\\
(iii) for any integer $N \geq 1,$ there exists a constant $ C_N>0$ such that 
$$\sup_{s \geq 0}\sup_{|x|\leq N}|\sigma_2(s,x)|\leq C_N.$$ 
Then the equation (1.1) has a unique strong solution.}
\vsp
\NI{\bf Remarks 1.1:} Suppose the conditions stated in Theorem 1.3 hold. Following Lemma 1 of Khlifa (2016) and arguments in Skorokhod (1965), it can be shown that the process $\{X_t,0\leq t \leq T\}$ satisfying the SDE (1.1) has the property (cf. Theorem 5.4, Klebaner (2012))
\be
E(\sup_{0\leq t \leq T}X_t^2)\leq C( 1+E(X_0^2))
\ee 
for some constant $C>0.$
\vsp
We will denote any positive constant by C through out the discussion in the following sections and the constant may differ from one step to another.
\vsp
\newsection{Preliminaries}
Let $(\Omega,\clf,\bar \clf, P)$ be a complete probability space with filtration $\bar\clf=\{{\clf}_t, t \geq 0\}$ satisfying the standard assumptions. We assume that all the processes under discussion are adapted to the filtration $\bar \clf.$ Consider the SDE
\be
dX_t=S(X_t)dt+\epsilon \;\sigma_1(t,X_t)\sigma_2(t,Y_t)\;dW_t, X_0=x_o\in R, 0\leq t \leq T 
\ee
where $W$ is the standard Wiener process and the function $S(.)$ is not known. We assume that sufficient conditions, as stated in Theorem 
1.3, hold on the functions $S(x), \sigma_1(t,x), \sigma_2(t,x)$ and the process $Y$ so  that the equation (2.1) has a unique strong solution $\{X_t, 0\leq t \leq T\}.$ We will consider the problem of nonparametric estimation of the function $S(.)$ based on the observation $\{X_t,0\leq t \leq T\}.$ A special case of this problem is studied in Kutoyants (1994) when the functions $\sigma_1(.,.)$ and $\sigma_2(.,.)$ are identically equal to one.
\vsp
Suppose $\{x_t, 0 \leq t \leq T\}$ is the solution of the differential equation
\be
\frac{dx_t}{dt}=S(x_t) , x_0 , 0 \leq t \leq T. 
\ee
\vsp
Suppose the following condition holds.\\
$ (A_1)(i)$ : There exists a constant $L > 0 $  such that $ \left|S(x) -S(y) \right| \leq L |x-y|, x,y \in R.$\\
$(A_1)(ii)$: There exists a constant $K>0$ such that
$$(\sigma_1(t,x))^2 \leq K(1+|x|^2), 0\leq t \leq T, x \in R;$$
$(A_1)(iii)$ The function $\sigma_2(t,x)$ is bounded for $0\leq t \leq T, x \in R.$
\vsp
It is clear that the condition $(A_1)(i)$ implies that there exists a constant $M>0$ such that
$$|S(x)|\leq M(1+|x|, x \in R.$$
\vsp
Since the function $x_t,$ satisfies the ordinary differential equation (2.2), it follows that
$$ |S(x_t)-S(x_s)| \leq L|x_t-x_s|=L|\int_t^sS(x_u)du|\leq L_1|t-s|, t,s \in R$$
for some constant $L_1>0.$
\vsp
\NI{\bf Lemma 2.1. } {\it Suppose the functions $S(x), \sigma_1(t,x)$ and $\sigma_2(t,x)$ satisfy the condition $(A_1)(i)-(iii)$.
Let  $ X_t$ and $x_t$ be  the solutions of the equation (2.1) and (2.2) respectively. Then, with probability one,
\be
(a)|X_t-x_t|< e^{Lt} \epsilon \sup_{0\leq v \leq t}|\int_0^v \sigma_1(s,X_s)\sigma_2(s,Y_s)dW_s|
\ee
and there exists a positive constant $C$ such that}
\be
(b)\sup_{0 \leq t \leq T} E (X_t-x_t)^2 \leq C e^{2LT} T \epsilon^2 .
\ee
\vsp
\NI{\bf Proof of (a) :} Let $u_t=|X_t-x_t| $. Then, by $(A_1)(i)$, we have
\bea
u_t & \leq & \int^t_0 \left| S(X_v)-S(x_v) \right| dv + \epsilon \;|\int_0^t \sigma_1(s,X_s)\sigma_2(s,Y_s)dW_s|\\\nonumber
& \leq & L \int^t_0 u_v dv + \epsilon \sup_{0\leq v \leq t}|\int_0^v \sigma_1(s,X_s)\sigma_2(s,Y_s)dW_s|.
\eea
Applying the Gronwell lemma (cf. Lemma 1.12, Kutoyants (1994), p.26),  it follows that
\be
u_t \leq  \epsilon \sup_{0\leq v \leq t}|\int_0^v \sigma_1(s,X_s)\sigma_2(s,Y_s)dW_s| e^{Lt}. \\
\ee
\vsp
\NI{\bf Proof of (b) :} Let
$$V(t)=\int_0^t \sigma_1(s,X_s)\sigma_2(s,Y_s)dW_s, 0\leq t \leq T.$$
Observe that 
$$E[V^2(t)]= \int_0^tE[\sigma_1(s,X_s)\sigma_2(s,Y_s)]^2ds $$
and the last term is finite by the condition $(A_1)$ and the Remark 1.1. 
The process $\{V(t),0\leq t \leq T\}$ is a martingale adapted to the filtration generated by the Wiener process and it follows that
$$E[\sup_{0\leq t \leq T}V^2(t)]\leq 4 E[V^2(T)]$$ 
(cf. Klebaner (2012), Theorem 7.31, p.203). 
\vsp
From (2.3), we have ,
\bea
E(X_t-x_t)^2 & \leq & e^{2Lt} \epsilon^2 E (\sup_{0\leq u\leq t}\int_0^u \sigma_1(s,X_s)\sigma_2(s,Y_s)dW_s)^2 \\\nonumber
&\leq & 4 e^{2Lt} \epsilon^2\int_0^t E(\sigma_1(s,X_s)\sigma_2(s,Y_s))^2 ds\\\nonumber
&=& e^{2Lt} C \epsilon^2\int_0^tE(\sigma_1(s,X_s))^2ds\;\;\mbox{(by the condition} (A_1)(iii))\\\nonumber
&=& e^{2Lt} C \epsilon^2\int_0^t E(1+ |X_s|^2)ds\;\;\mbox{(by the condition} (A_1)(ii)) \\\nonumber
&=& e^{2Lt} C t \epsilon^2\sup_{0\leq s \leq T}E(1+ |X_s|^2)\\\nonumber
& = & e^{2Lt}\epsilon^2 C_4 t E(1+X_0^2) \;\; \mbox{(by the Remark 1.1)}.\\\nonumber
\eea
for some positive constant C. Hence
\be
\sup_{0 \leq t \leq T} E (X_t-x_t)^2  \leq C \; e^{2LT} T\epsilon^2 . \\
\ee
\vsp
\newsection{ Main Results}
Let $\Theta_0(L)$ denote the class of all functions $S(.)$ satisfying the condition $(A_1)(i)$ with the same bound $L$. Let
$\Theta_k(L) $ denote the class of all functions $S(.)$ which are uniformly bounded by the same constant $C$
and which are $k$-times differentiable satisfying the condition
$$|S^{(k)}(x)-S^{(k)}(y)|\leq L|x-y|, x,y \in R$$
for some constant $L >0.$ Here $g^{(k)}(.)$ denotes the $k$-th derivative of $g(.)$ at $x$ for $k \geq 0.$ If $k=0,$ we interpret $g^{(0)}$ as $g.$
\vsp
Let $G(u)$ \ be a  bounded function  with finite support $[A,B]$ satisfying the condition
\vsp
\NI{$(A_2)$} $G(u) =0\;\; \mbox{for}\;\; u <A \;\;\mbox{and}\;\; u >B, \;\;\mbox{and}\;\;  \int^B_A G(u) du =1.$
\vsp
It is obvious that the following conditions are satisfied by the function $G(.):$
\begin{description}
\item{(i)} $ \int^\infty_{-\infty} G^2(u) du < \infty;$ and \\
\item{(ii)}$\int^\infty_{-\infty} u^{2 (k+1)} G^2(u) du < \infty, k \geq 0.$\\
\end{description}
We define a kernel type estimator of the trend $S_t=S(x_t)$ \ as
\be
\widehat{S}_t = \frac{1}{\varphi_\epsilon}\int^T_0 G \left(\frac{\tau-t}{\varphi_\epsilon} \right) d X_\tau
\ee
where the normalizing function  $ \varphi_\epsilon \rightarrow 0 $ with $\epsilon^2
\varphi^{-1}_\epsilon  \rightarrow 0 $ \ as \ $ \epsilon \rightarrow 0. $ Let $E_S(.)$ denote the expectation when the function $S(.)$ is the trend function.
\vsp
\NI{\bf Theorem 3.1:}  {\it Suppose that the trend function $S(.) \in \Theta_0(L)$ and  the function
$ \varphi_\epsilon \rightarrow 0$  such that $ \epsilon^2 \varphi^{-1}_\epsilon \longrightarrow  0 $  as $\epsilon
\rightarrow  0$. Further suppose that the condition $(A_1)$  holds. Then, for any $ 0 < c
\leq d < T ,$ the estimator $\widehat{S}_t$ is uniformly
consistent, that is,}
\be
\lim_{\epsilon \rightarrow 0} \sup_{S(x) \in \Theta_0(L)} \sup_{c\leq t \leq d } E_S ( |\widehat{S}_t - S (x_t)|^2)= 0.
\ee
\vsp
In addition to the conditions $(A_1)-(A_3),$  assume that
\vsp
\NI{$(A_3)$}$ \int^\infty_{-\infty} u^j G(u)   du = 0 \;\;\mbox{for}\;\; j=1,2,...k \;\;\mbox{and}\;\;\int_ {-\ity}^{\ity} |G(u)|u^{k+1}du < \ity.$
\vsp
\NI{\bf Theorem 3.2:} {\it Suppose that the function $ S(x) \in \Theta_{k+1}(L)$ and  $
\varphi_\epsilon = \epsilon^{\frac{2}{2k+3}}.$ Then, under the
conditions $(A_2)$ and $ (A_3),$}
\be
\limsup_{\epsilon \rightarrow 0} \sup_{S(x) \in \Theta_{k+1}(L)}\sup_{c \leq t \leq d} E_S (| \widehat{S}_t - S(x_t)|^2)
\epsilon^{\frac{-4(k+1)}{2k+3}} \ < \infty.
\ee
\vsp
\NI{\bf Theorem 3.3:} {\it Suppose that the function $S(x) \in \Theta_{k+1}(L)$ and $ \varphi_\epsilon= \epsilon^{\frac{2}{2k+3}}.$
Then, under the conditions $(A_1)-(A_3)$, the asymptotic distribution of
$$ \epsilon^{\frac{-2(k+1)}{2k+3}} (\widehat{S}_t - S(x_t)) $$
has the mean
\be
m =\frac{S^{(k+1)} (x_t)}{(k+1) !} \int^\infty_{-\infty} G(u) u^{k+1}\ du  
\ee
and  the limiting distribution is the same as that of the random variable
$$\epsilon^{-\frac{1}{2k+3}}\int_0^TG \left(\frac{\tau-t}{\varphi_\epsilon}\right))\sigma_1(\tau,X_\tau)\sigma_2(\tau,Y_\tau)dW_\tau $$
as $\epsilon \raro 0.$}
\vsp
\newsection{Proofs of Theorems}
\NI{\bf Proof of Theorem 3.1 :}
Following the equation (3.1), we have
\bea
\;\;\;\\\nonumber
E_S[(\widehat{S}_t -S_t(x))^2]  &=& E_S|\frac{1}{\varphi_\epsilon}\int^T_0 G \left(\frac{\tau-t}{\varphi_\epsilon} \right)dX_\tau-S(x_t)|^2\\\nonumber
&=& E_S \{ \frac{1}{\varphi_\epsilon}  \int^T_0 G \left(\frac{\tau-t}{\varphi_\epsilon} \right) \left(S(X_\tau) -S(x_\tau) \right)  d \tau \\ \nonumber
& &+ \frac{1}{\varphi_\epsilon}
 \int^T_0 G \left(\frac{\tau-t}{\varphi_\epsilon}\right) S_\tau(x) d \tau- S_t(x)\\\nonumber
 && + \frac{\epsilon}{\varphi_\epsilon} \int^T_0 G \left(\frac{\tau-t}{\varphi_\epsilon} \right)
 \sigma_1(\tau,X_\tau)\sigma_2(\tau,Y_\tau)d W_\tau\}^2\\ \nonumber
 & \leq  & 3 E_S \left[ \frac{1}{\varphi_\epsilon}  \int^T_0 G \left(\frac{\tau-t}{\varphi_\epsilon} \right) (S(X_\tau) -S(x_\tau)) d\tau\right]^2 \\ \nonumber
 & & + 3 E_S \left[\frac{1}{\varphi_\epsilon} \int^T_0 G \left(\frac{\tau-t}{\varphi_\epsilon} \right)S(x_\tau) d\tau -S(x_t) \right]^2 \\ \nonumber
 & & + 3\frac{\epsilon^2}{\varphi_{\epsilon^2}} \left[ \int^T_0 (G \left(\frac{\tau-t}{\varphi_\epsilon}\right))^2 E_S[\sigma_1^2(t,X_\tau)\sigma_2^2(t,Y_\tau)] d\tau \right]\\\nonumber
 &= & I_1+I_2+I_3 \;\;\mbox{(say).}\\\nonumber
\eea
Note that
\bea
\;\;\;\\\nonumber
I_1 &= &3 E_S \left[ \frac{1}{\varphi_\epsilon} \int^T_0 G
\left(\frac{\tau-t}{\varphi_\epsilon} \right) (S(X_\tau) -S(x_\tau))
d\tau \right]^2 \\\nonumber
&= & 3E_S  \left[ \int^\infty_{-\infty} G(u)
\left(S(X_{t+\varphi_\epsilon u}) -S(x_{t+\varphi_\epsilon u})
)\right) du\right]^2\\\nonumber
& \leq & 3 (B-A)\int^\infty_{-\infty} G^2(u) L^2 E \left(
X_{t+\varphi_\epsilon u}-x_{t+\varphi_\epsilon u} \right)^2 \ du
\;\;\mbox{(by using the condition $(A_1)$)}\\\nonumber
& \leq & 3(B-A) \int^\infty_{-\infty} G^2(u) \;\;L^2 \sup_{0 \leq t +
\varphi_\epsilon u \leq T}E \left(X_{t+\varphi_\epsilon u}
-x_{t+\varphi_\epsilon u}\right)^2 \ du \\\nonumber
& \leq & C \epsilon^2 \;\;\mbox{(by using (2.4))}\\\nonumber
\eea
for some positive constant $C.$Furthermore
\bea
\;\;\;\\ \nonumber
I_2 &= & 3E_S \left[ \frac{1}{\varphi_\epsilon} \int^T_0 G\left(
\frac{\tau-t}{\varphi_\epsilon}\right) S (x_\tau) d \tau - S (x_t)
\right]^2 \\ \nonumber
& \leq & 3 E_S \left[ \int^\infty_{-\infty} G(u)
\left(S(x_{t+\varphi_\epsilon u})-S(x_t) \right)  \ du \right]^2
\\ \nonumber
& \leq & 3 L^2 \left[ \int^\infty_{-\infty} G(u) u \varphi_\epsilon du
\right]^2\\ \nonumber
& \leq  & C \varphi^2_\epsilon  \int^\infty_{-\infty} G^2(u) u^2 \ du
\\ \nonumber
& \leq  & C \varphi^2_\epsilon  \;\;\mbox{(by $(A_2)(ii)$)}.\\ \nonumber
\eea
for some positive constant $C$. Furthermore  the last term  tends to zero as  $\epsilon \rightarrow 0.$
In addition
\bea
\;\;\;\\ \nonumber
I_3 &= & \frac{3 \epsilon^2}{\varphi_\epsilon^2} E_S \left(\int^T_0
G \left(\frac{\tau-t}{\varphi_\epsilon} \right) \sigma_1(\tau,X_\tau)\sigma_2(\tau,Y_\tau)dW_\tau
\right)^2 \\\nonumber
& \leq & \frac{3 \epsilon^2}{\varphi_\epsilon^2} [\left(\int^T_0
\{G \left(\frac{\tau-t}{\varphi_\epsilon} \right)\}^2 E(\sigma_1(\tau,X_\tau)\sigma_2(\tau,Y_\tau))^2d\tau\right)\\\nonumber
& \leq & \frac{6 \epsilon^2}{\varphi^2_\epsilon} C
 \left[\int^T_0 \left\{ G \left(\frac{\tau-t}{\varphi_\epsilon}
\right)\right\}^2 [E(\sigma_1(\tau,X_\tau))^2]d \tau  \right] \\\nonumber
&&\;\;\;\;\mbox{(from the boundedness of the function}\;\;\sigma_2(.,.))\\\nonumber
& \leq & \frac{C \epsilon^2}{\varphi^2_\epsilon} \varphi_\epsilon (\;\; \mbox{(by using $(A_1)(ii)$ and the Remark 1.4 )} \\ \nonumber
& = & C\frac{\epsilon^2}{\varphi_{\epsilon}} \\ \nonumber
\eea
for some positive constant $C$. Theorem 3.1 is now proved by using the equations (4.1) to (4.4).
\vsp
\NI {\bf Proof of Theorem 3.2 :}
By the Taylor's formula, for any $x \in R,$
$$ S(y) = S(x) +\sum^k_{j=1} S^{(j)} (x) \frac{(y-x)^j}{j !} +[ S^{(k)} (z)-S^{(k)} (x)] \frac{(y-x)^k}{k!} $$
for some $z$ such that $|z-x|\leq |y-x|.$
\vsp
Using this expansion, the equation (4.1) and the conditions  in the expression $I_2$ defined in the proof of  Theorem 3.1, it follows that
\ben
\;\;\\\nonumber
I_2 & \leq & 3 \left[
\int^\infty_{-\infty} G(u) \left(S_{t+\varphi_\epsilon u})(x) - S_t
(x) \right)  \ du \right]^2\\ \nonumber
&= & 3[ \sum^k_{j=1} S^{(j)}_t (x) (\int^\infty_{-\infty}
G(u) u^j du )\varphi^j_\epsilon ( j !)^{-1}\\\nonumber
& & \;\;\;\;+(\int^\infty_{-\infty}
G(u) u^k (S^{(k)}(z_u) -S^{(k)}_t (x))du
\varphi^k_\epsilon (k !)^{-1}]^2\\ \nonumber
\een
for some $z_u$ such that $|x_t-z_u|\leq |x_{t+\varphi_\epsilon u}-x_t| \leq C|\varphi_\epsilon u|.$ Hence
\bea
I_2 & \leq & CL^2 \left[  \int^\infty_{-\infty} |G(u)
u^{k+1}|\varphi^{k+1}_\epsilon (k!) ^{-1}  du  \right]^2
\\ \nonumber
& \leq & C(B-A)(k!)^{-2} \varphi^{2(k+1)}_\epsilon
\int^\infty_{-\infty} G^2(u) u^{2 (k+1)}\ du  \leq C\varphi_\epsilon^{2(k+1)}\\ \nonumber
\eea
for some positive constant $C$. Combining the relations (4.2) , (4.4) and
(4.5), we get that there exists a positive constant $C$ depending on
$$ \sup_{c \leq t \leq d}E_S
|\widehat{S}_t-S(x_t)|^2 \leq C (\epsilon^2\varphi^{-1}_\epsilon +  \varphi^{2(k+1)}_\epsilon +
\epsilon^2). $$ Choosing $ \varphi_\epsilon =
\epsilon^{\frac{2}{2k+3}},$  we get that $$ \limsup_{\epsilon
\rightarrow 0} \sup_{S(x) \in \Theta_{k+1} (L) } \sup_{c \leq t
\leq d} E_S |\widetilde{S}_t - S (x_t)|^2\epsilon^ {-
\frac{4(k+1)}{{2k+3}}} < \infty. $$ This completes the proof of
Theorem 3.2. \vsp \NI{\bf Remarks :} Choosing $\varphi_\epsilon =
\epsilon^{\frac{2}{3}}$  and without assuming the condition
$(A_3),$  it can be shown that $$ \limsup_{\epsilon \rightarrow 0}
\sup_{S(x) \in \Theta_0(L)} \sup_{c \leq t \leq d} E_S
|\widehat{S}_t - S(x_t)|^2 \epsilon^{-\frac{4 }{3}} < \infty $$
which gives a  slower  rate of convergence than the one obtained
in Theorem 3.2 .
\vsp \NI{\bf Proof of Theorem 3.3:}

From the equation (3.1), we obtain that

\bea
\;\;\;\;\\\nonumber
\lefteqn{ \widehat{S}_t -S(x_t)}\\\nonumber
 &= &[ \frac{1}{\varphi_\epsilon} \int^T_0 G \left(\frac{\tau-t}{\varphi_\epsilon} \right)
 \left( S(X_\tau)-S(x_\tau)\right) \  d \tau \\ \nonumber
 & & + \frac{1}{\varphi_\epsilon} \int^T_0 G \left( \frac{\tau-t}{\varphi_\epsilon}\right) S(x_\tau) d\tau -S(x_t)+ \frac{\epsilon}{\varphi_\epsilon} \int^T_0 G \left( \frac{\tau-t}{\varphi_\epsilon}\right) \sigma_1(t,X_\tau) \sigma_2(t,Y_\tau)dW_\tau]\\ \nonumber
 &= & \int^\infty_{-\infty} G(u) (S (X_{t+\varphi_\epsilon}u) - S (x_{t+\varphi_\epsilon
 u})) \ du  \\ \nonumber
 & & +\int^\infty_{-\infty} G(u) (S(x_{t+\varphi_\epsilon u})-
 S(x_t)) \ du \\ \nonumber
 & &+ \frac{\epsilon}{\varphi_{\epsilon}}\int^T_0 G
 \left(\frac{\tau-t}{\varphi_\epsilon}
 \right) \sigma_1(\tau,X_\tau) \sigma_2(\tau,Y_\tau) dW_\tau]\\\nonumber
&=& R_1+R_2+R_3 \;\;\mbox{(say).}\\\nonumber
\eea
Let $\varphi_\epsilon =\epsilon^{\frac{2}{2k+3}}$. By the Taylor's formula, for any $x \in R,$
$$ S(y) = S(x) +
\sum^{k+1}_{j=1} S^{(j)} (x) \frac{(y-x)^j}{j !} +
[ S^{(k+1)} (z)-S^{(k+1)} (x)
] \frac{(y-x)^{k+1}}{k+1!} $$
for some $z$ such that $|z-x|\leq |y-x|.$ Using this expansion, the equation (3.2) and the condition $(A_3),$ it follows that
\ben
(R_2-m)^2 & = & \left[\int^\infty_{-\infty} G(u) \left(S^{(k+1)} (z_u) - S^{(k+1)} (x_t)\right)\frac{(\varphi_\epsilon u)^{k+1}}{(k+1)!} \ du \right]^2
\\\nonumber
\een
for some $z_u$ such that $|x_t-z_u|\leq |x_{t+ \varphi_\epsilon u}-x_t| \leq C|\varphi_\epsilon u|$  and $m$ is as defined by the equation (3.4). Hence, by arguments similar to those given following the inequality (4.5), it follows that there exists a positive constant $C$ such that
\bea
(R_2-m)^2 & \leq & C L^2 \left(  \int^\infty_{\infty} G(u) u^{k+2} \frac{
\varphi^{k+2}_\epsilon}{(k+1)!} \ du \right)^2  \;\; (\mbox{by} (A_2))\\\nonumber
& \leq & C \varphi^{2(k+2)}_\epsilon  \;\;(\mbox{by} \;\;(A_3))\\\nonumber
\eea
for some positive constant $C$.Therefore
$$
\epsilon^{ - \frac{4(k+1)}{2k+3}}(R_2 -m)^2 \rightarrow 0 \;\;\mbox{as}\;\; \ \epsilon \rightarrow 0.
$$
Furthermore
\ben
0\leq \epsilon^{ - \frac{4(k+1)}{2k+3}}E_S[R^2_1] = \epsilon^{ - \frac{4(k+1)}{2k+3}} O(\phi_{\epsilon}^{2(k+2)})
\een
by arguments similar to those given for proving the inequality (4.5). Hence
$$
\epsilon^{ - \frac{4(k+1)}{2k+3}}E_S[R^2_1] \rightarrow 0\;\; \mbox{as}\;\; \epsilon \rightarrow
0.
$$
In addition, it follows that $ E_S [R^2_3]  $   is finite by $(A_2)(iii)$ and the
limiting distribution of the random variable  of the random variable
$$ \epsilon^{-\frac{1}{2k+3}} (\widehat{S}_t - S(x_t)) $$
is the same as that of the limiting distribution of the random variable
$$\epsilon^{-\frac{1}{2k+3}}\int_0^TG \left(\frac{\tau-t}{\varphi_\epsilon}\right))\sigma_1(\tau,X_\tau)\sigma_2(\tau,Y_\tau)dW_\tau $$
as $\epsilon \raro 0$ by an application of the Slutsky's lemma. This proves Theorem 3.3.
\vsp
\NI{\bf Funding :}\\ 

This work was supported under the scheme ``INSA Honorary Scientist" by the Indian National Science Academy  while the author was at the CR Rao Advanced Institute for Mathematics, Statistics and Computer Science, Hyderabad 500046, India.
\vsp
\NI{\bf References}
\begin{description}
\item Khlifa, M.B.H., Mishura,Y., Ralchenko, K., and Zili, M. (2016) drift parameter estimation in stochastic differential equation with multiplicative stochastic volatility, {\it Modern Stochastics: Theory and Methods}, {\bf 3}, 269-285.
\item Skorokhod, A.V. (1965) {\it Studies in the Theory of Random Processes}, Addison-Wesley Publishers, Reading, Mass.
\item Klebaner, F.C. (2012) {\it Introdcution to Stochatic Calculus with Applications}, Imperial College Press, London.
\item Kutoyants, Y. (1994) {\it Identification of Dynamical Systems with Small Noise}, Kluwer Academic Publishers, Dordrecht.

\end{description}
\end{document}